\documentclass{amsart}
\copyrightinfo{2006}{American Mathematical Society}
\usepackage{mathrsfs,amscd}
\usepackage{amssymb}
\usepackage{amssymb, amsmath, mathrsfs}
\usepackage{amscd}
\usepackage{latexsym}
\usepackage{amssymb}
\usepackage{amsmath}
\usepackage{amsthm}
\usepackage{indentfirst}

\newtheorem{theorem}{Theorem}[section]

\theoremstyle{definition}
\newtheorem{definition}[theorem]{Definition}

\theoremstyle{remark}
\newtheorem{remark}[theorem]{Remark}

\numberwithin{equation}{section}
\allowdisplaybreaks
\begin{document}

\title{Explicit self-similar solutions of time-like extremal hypersurfaces in $\mathbb{R}^{1+3}$}

\author[W. Yan]{Weiping Yan}
\address{School of Mathematics, Xiamen University, Xiamen 361000, People's Republic of China}
\email{yanwp@xmu.edu.cn}
\thanks{}

\subjclass[2010]{}
\date{June 2017}
\keywords{wave equation, self-similar solution}

\begin{abstract}
In this letter, two explicit self-similar solutions to a graph representation of time-like extremal hypersurfaces in Minkowski spacetime $\mathbb{R}^{1+3}$ are given. Meanwhile, there is an untable eigenvalue in the linearized time-like extremal hypersurfaces equation around two explicit self-similar solutions.
\end{abstract}
\maketitle

\section{Introduction and main results}
Let $\mathcal{M}$ be a timelike $(M+1)$-dimensional hypersurface, and $(\mathbb{R}^{D},g)$ be a $D$-dimensional Minkowski space, and $g$ be the Minkowski metric. At any time $t$, the spacetime volume in $\mathbb{R}^{D}$ of timelike hypersurface $\mathcal{M}$ can be described as a graph over $\mathbb{R}^M$, which satisfies
\begin{equation}\label{E1-0}
\mathcal{S}(u)=\int_{\mathbb{R}^M}\sqrt{1-|\partial_t u|^2+|\nabla u|^2}d^Mxdt.
\end{equation}
Critical points of action integral (\ref{E1-0}) give rise to submanifolds $\mathcal{M}\subset\mathbb{R}^{D}$ with
vanishing mean curvature, i.e. time-like extremal hypersurfaces. The Euler-Lagrange equation of (\ref{E1-0}) is
\begin{equation}\label{E1-1}
(1-u_{\alpha}u^{\alpha})\square u+u^{\beta}u^{\alpha}u_{\alpha\beta}=0,
\end{equation}
where $\forall\alpha,\beta=0,1,2,\ldots,M$, $u_{\alpha}=\frac{\partial u}{\partial x^{\alpha}}$, $z_{\alpha\beta}=\frac{\partial^2 u}{\partial x^{\alpha}\partial x^{\beta}}$ and $\square u=u_{\alpha\beta}g^{\alpha\beta}$.

Let $M=2$, $D=4$ and $r=|x|$. Then we obtain the radially symmetric membranes equation
\begin{equation}\label{E1-2}
u_{tt}-u_{rr}-\frac{u_r}{r}+u_{tt}u_r^2+u_{rr}u_t^2-2u_tu_ru_{tr}+\frac{1}{r}u_ru_t^2-\frac{1}{r}u_r^3=0.
\end{equation}
Above equation has been derived by Eggers and Hoppe \cite{Hop1}. Obviously, it is a quasilinear one dimensional wave equation. It is completely different with the case of wave map. There is a semilinear radially symmetric wave equation when reducing wave map (e.g. see \cite{Cos2}).
It is easy to see that (\ref{E1-2}) exhibits the following scaling invariance for any $\lambda>0$,
\begin{equation}\label{E1-3}
u(t,r)\mapsto u_{\lambda}(t,r)=\lambda u(\lambda^{-1}t,\lambda^{-1}r).
\end{equation}
Under this acaling the conserved energy
\begin{equation*}
E(u)=\int_{0}^{\infty}(\frac{1}{2}u_t^2+\frac{1}{2}u_r^2+F(u_t,u_r))rdr,
\end{equation*}
where $F'(u_t,u_r)=u_{tt}u_r^2+u_{rr}u_t^2-2u_tu_ru_{tr}+\frac{1}{r}u_ru_t^2-\frac{1}{r}u_r^3$, can be transformed as
\begin{equation*}
E(u_{\lambda})=\lambda E(u).
\end{equation*}
This means that the radially symmetric membranes equation is energy supercritical. So one expects smooth finite energy initial data to lead to finite time blow up, and the blow up rate is like self-similar blow up solution. One can compare this aspect with wave map, e.g. one can see \cite{B1,BB,Cos2,Cos3}. The investigation of blow up to membrane equation has attracted many work, one can see \cite{Hop1,Hop2,tian}.
Here we expect that the radially symmetric membranes equation (\ref{E1-2}) has self-similar blow up solutions. Eggers and Hoppes \cite{Hop1} gave a detail discussion on the existence of self-similar blow up solutions (not explicit self-similar solutions) to the radially symmetric membranes equation (\ref{E1-2}). Meanwhile, they gave some numerical analysis results on singularities of (\ref{E1-1}). In this note, we show two explicit self-similar blow up solutions to the radially symmetric membranes equation (\ref{E1-2}). After that, there is an untable eigenvalue in the linearized time-like extremal hypersurfaces equation around the explicit self-similar solutions. More precisely, our results are the following theorem.

\begin{theorem}
The radially symmetric membranes equation (\ref{E1-2}) has two explicit self-similar solutions
\begin{equation*}
\pm(T-t)\sqrt{1-(\frac{r}{T-t})^2}.
\end{equation*}
Furthermore, the linearized radially symmetric membranes equation around them has an unstable eigenvalue
$\nu=4$.
\end{theorem}

We remark that there is no explicit self-similar solution to the timelike $(1+1)$-dimensional hypersurface of $3$ dimensional Minkowski space, it can be described by the Lagrange equation
\begin{equation}\label{Y0}
u_{tt}-u_{xx}+u_{tt}u_x^2+u_{xx}u_t^2-2u_{t}u_{x}u_{tx}=0,
\end{equation}
which is an one dimensional quasilinear wave equation, and appeared in the paper of \cite{Hop1}.
It is so called Born-Infeld equation [1] describes the motion of a
string in the plane, and has been studied in \cite{BA1,BA2,BA3}.
In \cite{Hop1}, they made the self-similar ansatz
\begin{equation}\label{Y1}
u(t,x)=u_0-\hat{t}+\hat{t}^a\tilde{u}(\frac{x}{t^{b}})+\ldots,
\end{equation}
to above one dimensional quasilinear wave equation (\ref{Y0}), where $\hat{t}=t_0-t\rightarrow0$ when $t\rightarrow t_0$, $u_0=u(0,x)$, $a$, $b$ and $t_0$ are fixed constant.
Then a self-similar solution like (\ref{Y1}) to (\ref{Y0}) has been shown. Furthermore, they also showed that there is self-similar singularity like (\ref{Y1}) to the timelike $(1+2)$-dimensional hypersurface of $4$ dimensional Minkowski space, i.e. the radially symmetric membranes equation (\ref{E1-2}). In this note, two explicit self-similar solutions to equation (\ref{E1-2}) are obtained.

\section{Two explicit self-similar solutions and an untable eigenvalue}
\subsection{Two explicit self-similar solutions}
Self-similar solutions are invariant under the scaling (\ref{E1-3}), so let
\begin{eqnarray*}
&&\rho=\frac{r}{T-t},\\
&&u(t,r)=(T-t)\phi(\rho),
\end{eqnarray*}
where $T$ is a positive constant, which is a parameter. Inserting this ansatz into equation (\ref{E1-2}), one gets a quasilinear ordinary differential equation
\begin{equation}\label{E2-1}
\rho(1-\rho^2)\phi''+\phi'-\phi'\phi^2+2\rho\phi(\phi')^2-\rho\phi''\phi^2+(1-\rho^2)(\phi')^3=0.
\end{equation}
We are only interested in smooth solutions in the backward light-cone of the blow up point $(t,r)=(T,0)$, i.e. $\rho$ in the closed interval $[0,1]$. This solution is such that
\begin{equation}\label{E2-1R1}
\frac{\partial^n}{\partial r^n}\phi(\frac{r}{T-t})|_{r=0}=(T-t)^{-n}\frac{d^n\phi}{dy^n}(0),
\end{equation}
where when the $n$-th derivative is even, there is $\frac{d^n\phi}{dy^n}(0)\neq0$, this means it is diverges as $t\rightarrow T$. When $n$ is odd, there is $\frac{d^n\phi}{dy^n}(0)=0$. This condition is different with wave map (e.g. see \cite{BB}).
One can see each self-similar solution $\phi(\rho)\in\mathbb{C}^{\infty}[0,1]$ describes a singularity developing in finite time from smooth initial data. Since (\ref{E2-1}) is quasilinear ODE, there are some difficulties to prove the existence of smooth of solutions.
But from the structure of nonlinear term in equation (\ref{E2-1}), it can be rewritten as
\begin{equation}\label{E2-1R}
\rho(1-\rho^2-\phi^2)\phi''+\phi'-\phi'\phi^2+2\rho\phi(\phi')^2+(1-\rho^2)(\phi')^3=0.
\end{equation}
We found that if we let
\begin{equation*}
1-\rho^2-\phi^2=0,
\end{equation*}
then there holds
\begin{equation*}
\phi'-\phi'\phi^2+2\rho\phi(\phi')^2+(1-\rho^2)(\phi')^3=0.
\end{equation*}
Hence we obtain two explicit solutions
\begin{equation}\label{E2-2}
\phi(\rho)=\pm\sqrt{1-\rho^2},
\end{equation}
which satisfies condition (\ref{E2-1R1}).

Consequently, two explicit self-similar solutions of equation (\ref{E1-2}) are
\begin{equation}\label{E2-3}
\phi^T(t,r)=\pm(T-t)\sqrt{1-(\frac{r}{T-t})^2},
\end{equation}
which exhibit smooth for all $0<t<T$, but which break down at $t=T$ in the sense that
\begin{eqnarray*}
\partial_{rr}\phi^T(t,r)|_{r=0}&=&\pm((T-t)^2-r^2)^{-\frac{1}{2}}|_{r=0}\pm r^2((T-t)^2-r^2)^{-\frac{3}{2}}|_{r=0}\nonumber\\
&=&\pm\frac{1}{T-t}\rightarrow+\infty,~~as~~t\rightarrow T^{-},
\end{eqnarray*}
and the dynamical behavior of them are as attractors.

On the other hand, from the form of $\phi^T(t,r)$ in (\ref{E2-3}), it requires that
\begin{equation*}
1-(\frac{r}{T-t})^2\geq0.
\end{equation*}
So we consider the dynamical behavior of self-similar solutions in
of the backward lightcone
\begin{equation*}
\mathcal{B}_T:=\{(t,r):t\in(0,T),~~r\in[0,T-t]\}.
\end{equation*}

\begin{remark}
Obviouly, self-similar solutions (\ref{E2-3}) are cycloids. The sphere begins to expand until it starts to shrink and eventually collapses to a point in a finite time $\tilde{T}=T-r_0$, i.e. $\phi^{T}(\tilde{T},r_0)=0$.
Here $r_0$ is a fixed positive constant in the backward lightcone.
\end{remark}

\subsection{Mode instability of self-similar solutions $\phi^T(t,r)$}

We introduce the similarity coordinates
\begin{equation}\label{E2-5}
\tau=-\log(T-t),~~\rho=\frac{r}{T-t},
\end{equation}
then we denote by
\begin{equation*}
\tilde{v}(\tau,\rho)=e^{\tau}u(T-e^{-\tau},\rho e^{-\tau}),
\end{equation*}
equation (\ref{E1-2}) is transformed into
\begin{eqnarray}\label{E2-4}
\tilde{v}_{\tau\tau}-\tilde{v}_{\tau}&-&(1-\rho^2)\tilde{v}_{\rho\rho}-\frac{1}{\rho}\tilde{v}_{\rho}+2\rho \tilde{v}_{\tau\rho}+\tilde{v}_{\rho}^2(\tilde{v}_{\tau\tau}+\tilde{v}_{\tau}-2\tilde{v})+\tilde{v}_{\rho\rho}
(\tilde{v}-\tilde{v}_{\tau})^2\nonumber\\
&&-2\tilde{v}_{\rho}\tilde{v}_{\tau\rho}(\tilde{v}_{\tau}-\tilde{v})+\frac{1}{\rho}\tilde{v}_{\rho}(\tilde{v}_{\tau}-\tilde{v})^2+\frac{1}{\rho}(\rho^2-1)\tilde{v}_{\rho}^3=0.
\end{eqnarray}
In the similarity coordinates (\ref{E2-5}), the blow up time $T$ is changed to $\infty$. So the stability of blow up solutions of quasilinear wave equation (\ref{E1-2}) as $t\rightarrow T^{-}$ is transformed into the asymptotic stability of quasilinear equation (\ref{E2-4}) as $\tau\rightarrow\infty$. We consider the problem of mode stability or instability of linear equation of (\ref{E2-4}). Let solution of (\ref{E2-4}) takes the form
\begin{eqnarray}\label{E2-4R}
\tilde{v}(\tau,\rho)&=&\phi(\rho)+v(\tau,\rho),
\end{eqnarray}
where $\phi(\rho)$ is defined in (\ref{E2-2}).

Inserting this ansatz (\ref{E2-4R}) into equation (\ref{E2-4}), $v(\tau,\rho)$ satisfies one dimensional equation
\begin{equation}\label{E2-6r}
\begin{array}{lll}
(1+\phi'^2)v_{\tau\tau}&-&(1-\phi'^2+2\phi''\phi+\frac{2}{\rho}\phi')v_{\tau}+2(\phi\phi'+\rho)v_{\tau\rho}-(1-\rho^2-\phi^2)v_{\rho\rho}\\
&&-\frac{1}{\rho}(1+4\rho\phi'\phi-3(\rho^2-1)\phi'^2-\phi^2)v_{\rho}+\frac{1}{\rho}(-2\rho\phi'^2+2\rho\phi''\phi+2\phi'\phi)v\\
&&-2\phi v_{\rho}^2+v_{\rho}(v_{\tau\tau}+v_{\tau}-2v)(v_{\rho}+2\phi')+\phi''(v-v_{\tau})^2\\
&&+v_{\rho\rho}[(v-v_{\tau})^2+2\phi(v-v_{\tau})]-2\phi'v_{\tau\rho}(-v+v_{\tau})+2\phi v_{\rho}v_{\rho\tau}\\
&&-2v_{\rho}v_{\rho\tau}(-v+v_{\tau})+\frac{1}{\rho}v_{\rho}[(v-v_{\tau})^2+2\phi(v-v_{\tau})]+\frac{1}{\rho}\phi'(v-v_{\tau})^2\\
&&+(\rho-\frac{1}{\rho})(v_{\rho}^3+3\phi'v_{\rho}^2)=0.
\end{array}
\end{equation}
From the exact form of $\phi(\rho)$ in (\ref{E2-2}), it has
\begin{eqnarray*}
&&1-\rho^2-\phi^2=0,\\
&&\phi\phi'+\rho=0.
\end{eqnarray*}
Thus one dimensional equation (\ref{E2-6r}) is loss of hyperbolicity, the linear equation takes the form
\begin{equation}\label{E2-7r}
\begin{array}{lll}
v_{\tau\tau}&+&3v_{\tau}-4v=0.
\end{array}
\end{equation}
then let $v(\tau,\rho)=e^{\nu\tau}u_{\nu}$, it leads to a eigenvalue problem
\begin{equation}\label{E2-4R1}
(\nu^2+3\nu-4)u_{\nu}=0.
\end{equation}
As in \cite{Cos2, Cos3}, we introduce defnitions of mode stable and untable of solution $u_{\nu}$ to equation (\ref{E2-4R1}).

\begin{definition}
A non-zero smooth solution $u_{\nu}$ of (\ref{E2-4R1}) is called mode stable if $Re~\nu<0$ holds. The eigenvalue $\nu$ is called a stable eigenvalue. Otherwise, if $Re~\nu\geq0$ is called mode untable of non-zero smooth solution $u_{\nu}$ of (\ref{E2-4R1}). $\nu$ is called an unstable eigenvalue.
\end{definition}

From (\ref{E2-4R1}), it is easy to check that $\nu=4$ and $\nu=-1$ are two eigenvalues of linear operator to $(\ref{E2-7r})$.
By the definition of mode unstable, we know that two explicit self-similar solutions $\phi^T(t,r)$ in (\ref{E2-3}) of time-like extremal hypersurfaces equation (\ref{E1-2}) are mode unstable.

Now we propose a \textbf{question:}

Can we find a stability decomposition of blow up solutions to the radially symmetric membranes equation (\ref{E1-2})?

In the paper \cite{Yan}, we will answer this question, and show that there is a stability decomposition of blow up solutions concerning with two self-similar solutions found in Theorem 1.1.

\textbf{Acknowledgments.} The author expresses his sincere thanks to Prof. G. Tian for his many kind helps and suggestions. The author also expresses his sincerely
thanks to Prof. J. Hoppe for his pointing out two explicit solutions being lightlike, and his suggestion and informing me his interesting papers \cite{hop}.

\end{document}